\documentclass[leqno,a4paper,11pt]{article}
\usepackage{amsmath}
\begin{document}
\newtheorem{izrek}{Theorem}[section]
\newtheorem{posledica}[izrek]{Corollary}
\newtheorem{lema}[izrek]{Lemma}
\newtheorem{trditev}[izrek]{Proposition}
\newtheorem{korak}[izrek]{Step}
\newtheorem{rem}[izrek]{Remark}
\newtheorem{rems}[izrek]{Remarks}
\newcommand{\iy}{\infty}

\title{NONSURJECTIVE NEARISOMETRIES OF BANACH SPACES}

\author{Peter \v Semrl%
\footnote{The first author was supported in part by a grant from the
Ministry of
Science of Slovenia.}
\ and Jussi V\"ais\"al\"a}
\date{30.12.01}
\small
This work has been submitted to Academic Press for possible publication.
Copyright may be transferred without notice, after which this version may no
longer be accessible.
\normalsize
\maketitle

\begin{abstract}
We obtain sharp approximation results for into near\-isometries between
$L^p$ spaces and near\-isometries  into a Hilbert space.
Our main theorem is the optimal
approximation result for near\-surjective near\-isometries 
between
general Banach spaces.\\

\end{abstract}

\vskip 18pt
\section{Introduction}
Let $X$ and $Y$ be real Banach spaces and $\varepsilon\ge 0$. A map
$f: X \to Y$ is called an $\varepsilon$-{\it near\-isometry} if
$$
| \, \| f (x) - f(y) \| - \| x - y \| \, | \le \varepsilon
$$
for all $x,y \in X$. We also say that $f$ is a {\it near\-isometry} if it 
is an
$\varepsilon$-near\-isometry for some $\varepsilon \ge 0$.

The basic question is how close is
$f$ to an actual iso\-metry. Note that when studying such
maps there is no loss of generality in assuming that $f(0) = 0$. Indeed, if 
a mapping $f$
is an $\varepsilon$-near\-isometry then $f - f(0)$ is also an
$\varepsilon$-near\-isometry, and
$f - f(0)$ can be approximated by an iso\-metry $U$ if and only if $f$ is 
close to the iso\-metry $U + f(0)$.

The study of this problem was initiated by D.H. Hyers and S.M. Ulam \cite{HyU},
who proved that for every surjective $\varepsilon$-near\-isometry $f:X\to 
Y$ between
real Hilbert spaces
satisfying $f(0) = 0$
there is a bijective linear iso\-metry $U:X \to Y$
such that
\begin{equation}\label{prvae}
\| f(x) - Ux \| \le 10 \varepsilon
\end{equation}
for every $x\in X$. In 1983, after many partial results extending over 
almost four decades,
J. Gevirtz \cite{Gev} extended this theorem to arbitrary Banach spaces $X$ 
and $Y$ with the
better estimate $5\varepsilon$ in (\ref{prvae}). Finally, M. Omladi\v c and 
P. \v Semrl
\cite{OmS} showed that $2\varepsilon$ is a sharp constant in (\ref{prvae}) 
for general Banach
spaces. The surjectivity assumption is indispensable in this result. 
Namely, already
Hyers and Ulam \cite{HyU} gave an example of a near\-isometry $f$ from the 
real line into the
Euclidean plane such that $\{ \| f(x) - Ux \|\, : \, x\in {\bf R} \}$ is an 
unbounded set
for every iso\-metry $U : {\bf R}\to {\bf R}^2$. Another such example is 
given by
$f(x) = (x, \sqrt{2\varepsilon |x|})$.

The classical Mazur-Ulam theorem \cite{MaU} states that every onto 
iso\-metry $f: X \to
Y$ with $f(0) = 0$ is linear. When studying into near\-isometries we have 
to start with
a nonsurjective substitute for the Mazur-Ulam theorem. Such a substitute 
was given by
T. Figiel (see \cite{Fig} or \cite[14.2]{BeL}), who proved that for any 
iso\-metry $f:
X
\to Y$ with
$f(0) = 0$ there is a linear operator $T$ of norm one from $\overline{ {\rm 
span}}\,
f(X)$ onto
$X$ such that $T\circ f$ is the identity on $X$. S. Qian \cite{Qia}
gave an example to show that the approximate version of Figiel's theorem 
does not hold
for general Banach spaces. He also proved the approximate version of 
Figiel's theorem in two
special cases when $X$ and $Y$ are $L^p$ spaces or when $Y$ is a Hilbert 
space. More precisely,
in these two special cases the assumption that
$f:X \to Y$ is an into $\varepsilon$-near\-isometry with $f(0) = 0$
yields the existence of a continuous linear operator $T: Y \to X$ with $\| 
T \| = 1$ such that
\begin{equation}\label{SQ}
\| Tf(x) - x \| \le 6\varepsilon
\end{equation}
for every $x\in X$. In the next section we shall obtain a sharp version of 
these two theorems
by reducing $6\varepsilon$ in (\ref{SQ}) to $2\varepsilon$.

S.J. Dilworth \cite{Dil} noted that a near\-isometry
$f:X \to Y$
which can be uniformly approximated by a linear iso\-metry
maps $X$ approximately onto some closed subspace of $Y$. More precisely,
let $f:X\to Y$ be a map, $Y_1$ a closed subspace of $Y$,
and $\delta$ a nonnegative real number. Then we say that $f$ maps $X$
$\delta$-{\it onto} $Y_1$ if for every $y\in Y_1$ there is $x\in X$ with
$\| f(x) - y \| \le \delta$ and for every $u\in X$ there is $v\in Y_1$
with $\| f(u) -v \| \le \delta$.
When
$f$ maps $X$ $\delta$-onto $Y$  we simply say that 
$f$ is
$\delta$-{\it near\-surjective}. Moreover, $f$ is {\it near\-surjective} if 
it is
$\delta$-near\-surjective for some $\delta \ge 0$. Now, if
$f:X
\to Y$ is a near\-isometry that can be uniformly approximated by a linear 
iso\-metry
$U$ (that is, there exists a nonnegative real number $\delta$ such that $\| 
f(x) -
Ux
\| \le
\delta$,
$x\in X$) then clearly $f$ maps $X$ $\delta$-onto
the closed subspace $U(X)\subset Y$. The converse is also true.
Namely, Dilworth \cite[Prop. 2]{Dil} proved that if
$\varepsilon, \delta \ge 0$ and
$f:X \to Y$ is an $\varepsilon$-near\-isometry between
Banach spaces satisfying $f(0)=0$ and mapping
$X$ $\delta$-onto some closed subspace $Y_1\subset Y$, then
there exists a bijective linear iso\-metry $U: X \to Y_1$ such that
\begin{equation}\label{SD}
\| f(x) - Ux \| \le 12\varepsilon + 5\delta, \ \ \ x\in X.
\end{equation}
The natural question here is how far the estimate
$12\varepsilon + 5\delta$ is from the optimal one. If $\delta = 0$ then
we already know that $2\varepsilon$ is the sharp estimate in (\ref{SD}). The
example given by Dilworth \cite[Remark 1]{Dil} shows that the optimal estimate
in (\ref{SD}) increases at least linearly with $\delta$. So one may ask if the
estimate $12\varepsilon + 5\delta$ can be replaced by
$2\varepsilon + c\delta$ for some positive
$c$? Using the result of Dilworth and some geometrical considerations we
shall show
that the answer is affirmative
and that for general Banach spaces the optimal value of $c$ is 2. In the
case that
$Y$ is a Hilbert space the sharp bound in (3) is given by $2\varepsilon +
\delta$.

In our main result we consider nearsurjective nearisometries.
A lot of work has been done on surjective
near\-isometries. The notion of near\-isometries was motivated by the fact 
that
real-world observations have always some small error. So, the map $f$ which
associates to the real-world points the points in a
mathematical model is always a near\-isometry. But if we cannot measure the 
distances
exactly then we cannot
check if $f$ is surjective. Hence, we believe that it is more natural to
study the maps that are both near\-isometric and near\-surjective.
We shall prove that for every pair of Banach spaces $X$ and $Y$ and for
every $\delta$-near\-surjective $\varepsilon$-near\-isometry $f:X \to Y$ 
with $f(0) = 0$
there exists a bijective linear iso\-metry $U:X \to Y$ such that
\begin{equation}\label{twoplustwo}
\| f(x) - Ux \| \le 2\varepsilon , \ \ \ x\in X.
\end{equation}
Thus, the sharp estimate (\ref{twoplustwo}) is independent of $\delta$.
Up to now the best known approximation estimates were $12\varepsilon +
5\delta$ obtained by Dilworth \cite[Prop. 2]{Dil} and $2\varepsilon +
35 \delta$ obtained by J. Tabor \cite{Tab}. The statement of our main
result has also been known to be true under the additional condition
that the set of points in $X$ at which the norm is Fr\' echet differentiable
is dense \cite[Th. 2]{Dil} as well as under the stronger assumption that $f$
is surjective \cite{OmS}.\\

\section{Near\-isometries between $L^p$ spaces}

The definition of a near\-isometry obviously makes sense for maps between 
any metric
spaces. We start this section with a statement on near\-isometries defined 
on rays.
This result is of its own interest but will also be needed for
a sharp approximate version of Figiel's theorem for $L^p$ spaces.
\begin{izrek}\label{essay}
Suppose that $Y$ is a Banach space and that $f:[0, \infty) \to Y$ is an
$\varepsilon$-near\-isometry with $f(0)=0$. Then there is a linear functional
$F:Y \to {\bf R}$ of norm one such that $t- 2\varepsilon \le F f(t) \le t +
\varepsilon$ for all $t\ge 0$.
\end{izrek}
{\bf Proof.} We modify an idea of Figiel \cite{Fig}. For each integer $n>
\varepsilon$ we have $\| f(n) \| \ge n-\varepsilon > 0$. By the Hahn-Banach
theorem there is a linear functional $F_n : Y \to {\bf R}$ such that $\| 
F_n \| =1$
and $F_n f(n) = \| f(n)\| \in [n-\varepsilon , n+\varepsilon]$. For all
$t\ge 0$ we have $|F_n f(t)| \le \|f(t)\| \le t+\varepsilon$.

Let $0\le t \le n$. Since $\| F_n \| = 1$ and since $f$ is an 
$\varepsilon$-near\-isometry,
we obtain
$$
| F_n f(t) - F_n f(n) | \le \|f(t) - f(n)\| \le n-t+\varepsilon.
$$
Hence
\begin{align*}
F_n f(t) & = F_n f(n) - F_n (f(n) - f(t)) \\
& \ge n-\varepsilon -
(n-t+\varepsilon) = t-2\varepsilon.
\end{align*}
Since $\| F_n \| = 1$ for all $n$, it follows by Alaoglu's theorem that the 
sequence
$(F_n )$ has a cluster point $F\in Y^*$ in the $w^*$ topology. Then $\| F\| 
\le 1$ and
$$
t - 2\varepsilon \le \| F\|\| f(t)\| \le \| F\| (t + \varepsilon).
$$
Dividing by $t$ and letting $t \to\iy$ yields $\| F\| \ge 1$, and the 
theorem is proved.\\

Let $X$ be a uniformly convex Banach space with modulus of convexity $\tau$
(see \cite[p. 409]{BeL}).
This means that
$\tau : [0,2] \to [0,1]$ is an increasing function such that $\tau (s) >0$
for $0 < s \le 2$ and $\| x + y \| /2 \le 1 - \tau (s)$ whenever $\| x \| 
\le 1$,
$\| y \| \le 1$, and $\| x - y \| \ge s$. For our purposes, it is more 
convenient to
use the function $\gamma : [0,1] \to [0,2]$ defined by
\begin{equation}\label{gamma}
\gamma (t) = \sup \{ s\in [0,2] \, : \, \tau (s) \le t \}.
\end{equation}
We have $\gamma (t) \rightarrow 0$ as $t\rightarrow 0$. Indeed, let $0 < 
\varepsilon
\le 2$. For $s\ge \varepsilon$ we have $\tau (s) \ge \tau (\varepsilon )$.
Hence $\gamma (t) \le \varepsilon$ whenever $ 0 \le t < \tau (\varepsilon )$.

One can show that the minimal modulus of convexity $\tau_X$ \cite[p. 
409]{BeL} is a
homeomorphism onto $[0,1]$ and then $\gamma$ is simply $\tau_X^{-1}$, but 
we do not need
this fact.

For $R\ge 0$ and $x\in X$ set $\overline{B} (x, R) = \{ y\in X \, : \, \| 
y-x \|
\le R \}$. In the case $x=0$ we write shortly $\overline{B} (0, R)=
\overline{B} (R)$.
\begin{lema}\label{RRR}
Let $X$ be a uniformly convex Banach space with modulus of convexity $\tau$
and let $\gamma$ be as above. Assume that $0< r < R$ and that $x,y \in
\overline{B} (R)$ with $\| x + y \| / 2 \ge r$. Then $\| x-y \| \le
R\gamma ( 1 - r/R )$.
\end{lema}
{\bf Proof.} Suppose first that $R=1$. Assume that $\| x - y \| = s >
\gamma (1- r)$. Then $r\le \| x+y \| /2 \le 1 - \tau (s)$, and hence
$\tau (s) \le 1 - r$, which yields the contradiction
$\gamma ( 1-r) \ge s$.

In the general case we set $a=x/R$, $b=y/R$. Then $\| a+b \| /2 \ge r/R$.
By the special case we obtain $\| a - b \| \le \gamma ( 1- r/ R)$, and
the lemma follows.\\
\begin{trditev}\label{hus}
Suppose that $X$ and $Y$ are Banach spaces with $Y$ uniformly convex. Let
$f : X \to Y$ be a near\-isometry with $f(0) = 0$. Then the limit
$\varphi (x) = \lim_{ s \rightarrow \infty} f(sx)/s$ exists for each $x\in X$,
and the map $\varphi : X \to Y$ is a linear iso\-metry.
\end{trditev}
{\bf Proof.} Let $x$ be any vector in $X$. We want to prove that
 $\lim_{ s \rightarrow \infty} f(sx)/s$ exists. We may assume that
$\| x \| = 1$. Suppose that $f$ is an $\varepsilon$-near\-isometry.
Let $2\varepsilon < s < t$. Setting $y= sf(tx) / t$ we have
$$
\| y \| \le s (t+\varepsilon ) / t \le s + \varepsilon
$$
and
$$
\| y - f(tx) \| = | s / t - 1 | \, \| f(tx) \| \le (t+
\varepsilon ) (t-s) /t \le t-s + \varepsilon.
$$
Hence $y$ and $f(sx)$ lie in the convex set
$$
C=\overline{B} ( s+ \varepsilon ) \cap
\overline{B} ( f(tx), t-s + \varepsilon ).
$$
For each $z\in C$ we have
$$
\| z \| \ge \| f(tx) \| - \| z - f(tx) \| \ge
t - \varepsilon - ( t-s + \varepsilon) = s - 2 \varepsilon,
$$
and consequently, $C \subset
\overline{B} ( s+ \varepsilon ) \setminus B( s - 2\varepsilon )$. By Lemma 
\ref{RRR}
this implies that
$$
\| y - f(sx) \| \le ( s+ \varepsilon ) \gamma \left( 1 -
\frac{s - 2\varepsilon}{ s + \varepsilon} \right) =
(s+ \varepsilon ) \gamma \left( \frac{3\varepsilon }{s + \varepsilon} \right).
$$
Hence
$$
\left\| \frac{f(sx)}{ s} - \frac{f(tx)}{ t }\right\| = \| f(sx) -y \| / s \le
(1+ \varepsilon / s ) \gamma \left ( \frac{3\varepsilon }{ s + \varepsilon} 
\right)
\rightarrow 0
$$
as $s \rightarrow \infty$. Since $Y$ is complete, this implies the 
existence of
$\varphi (x)$.

Dividing the inequality
$$
| \, \| f (sx) - f(sy) \| - \| sx - sy \| \, | \le \varepsilon
$$
by $s$ and sending $s$ to infinity we see that $\varphi$ is an iso\-metry.
Clearly $\varphi (0) = 0$, and since $Y$ is strictly convex, $\varphi$
must be linear. This completes the proof.\\

Note that the above statement has been proved by D.G. Bourgin \cite[Theorem 
4]{Bou}
under two additional assumptions. We are now ready to prove the main result 
of this
section.
\begin{izrek}\label{FigielLp}
Let
$(\Omega_i , \Sigma_i , \mu_i )$, $i=1,2$, be two measure spaces, let $1 < 
p < \iy$, and
let $X_i = L^p (\Omega_i , \Sigma_i , \mu_i )$. Assume that $f: X_1 \to X_2$
is an $\varepsilon$-near\-isometry
with $f(0)=0$. Then there is a continuous linear operator
$T: X_2\to X_1$ with
$\| T \| = 1$ such that
$$
\| Tf(x)-x \| \le 2\varepsilon, \ \ \ x\in X_1.
$$
\end{izrek}
{\bf Proof.}
We shall several times make use of the fact that the norm function $x 
\mapsto \| x\|$ of
$X_2$ has a Fr\' echet derivative $j(z)$ at every nonzero point $z \in 
X_2$. Alternatively,
$j(z) \in X_2^*$ is the unique supporting functional at $z$, characterized 
by the
conditions $\| j(z)\| \le 1$ and $j(z)z = \|z\|$.

Suppose that $w \in X_2$ with $j(z)w = 0$. From the definition of Fr\' echet
diffe\-rentiability it easily follows that
\begin{equation}\label{A}
\lim_{t \to\iy} (\|tz+w\| - t\|z\|) = 0;
\end{equation}
see \cite[(2), p. 472]{Dil}.

Since $X_2$ is uniformly convex, Proposition \ref{hus} gives a linear 
iso\-metry $\varphi:
X_1 \to X_2$. Then $K = \varphi
(X_1)$ is a closed linear subspace of $X_2$, and $K$ is linearly isometric 
to $X_1$. We
shall consider $\varphi$ as a bijective linear iso\-metry of $X_1$ onto 
$K$. Then $\psi =
f\varphi^{-1}: K \to X_2$ is an $\varepsilon$-near\-isometry with
$\psi(0)=0$.

It
follows from \cite[p. 162]{Lac} that there is a linear projection $P$ from 
$X_2$ onto
$K$ with $\| P\| = 1$.
We shall show that $T = \varphi^{-1}P: X_2 \to X_1$ is the desired map. 
Clearly $\| T\| =
\| P\| = 1$. Let $x \in X_1$. Setting $y = \varphi x$ we have
$$
\| Tf(x) - x\| = \| \varphi^{-1}Pf(x) - x\| = \| Pf(x) - \varphi x\| = \| 
P\psi (y) - y\|.
$$
Consequently, it suffices to show that
$$
\| P\psi (y) - y\| \le 2\varepsilon.
$$

We first prove the weaker estimate
\begin{equation}\label{B}
\| P\psi (y) - y\| \le 3\varepsilon.
\end{equation}
Set $\lambda = \| P\psi (y) - y\|$ and choose a unit vector $z \in K$ with 
$y - P\psi (y) =
\lambda z$. By Theorem \ref{essay} there is a linear functional $F: X_2 \to 
{\bf R}$ with
$\| F\| = 1$ such that $|F\psi (tz) - t| \le 2\varepsilon$ for all $t \ge 
0$. For $t > 0$
we have
$$
|F(f(t\varphi^{-1}z)/t) - 1| = |F\psi (tz) - t|/t \le 2\varepsilon/t.
$$
As $t \to\iy$, this yields $Fz = 1$. Since $\| F\| = 1$, we have $F = 
j(z)$. Moreover,
since $\| FP\| \le 1$ and $FPz = Fz = 1$, we have $F=FP$.

We can write $P\psi (y) = \alpha z + w$ with $\alpha = F\psi(y)$ and $Fw = 
0$. Then $y
= (\lambda+\alpha)z+w$. Let
$t > 0$. Since $F\psi (tz) \ge t - 2\varepsilon$, we obtain
\begin{align*}
t - \alpha - 2\varepsilon &\le F\psi (tz) - F\psi (y) \le \| \psi (tz) - 
\psi (y)\| \le
\| tz - y\| +
\varepsilon\\
& = \| (t-\lambda-\alpha)z - w\|
+ \varepsilon.
\end{align*}
As $t \to\iy$, this and (\ref{A}) yield (\ref{B}).
We next improve (\ref{B}) to the sharp estimate $\lambda \le 2\varepsilon$. 
For each
positive integer $m$ we write
\begin{equation}\label{C}
\psi (y+mz) - \psi (y) = k_mz + w_m,
\end{equation}
where $k_m \in {\bf R}$ and $Fw_m = 0$. Setting $p_m = P\psi (y+mz) - 
(y+mz)$, applying
$P$ to (\ref{C}) and observing that $Pz=z$ we obtain
$$
k_mz + Pw_m = P\psi (y+mz) - P\psi (y) = p_m + mz + \lambda z.
$$
Applying $F$ to this equality and observing that $Fz = 1 $ and $FPw_m = 
Fw_m = 0$ we
get
$
k_m = Fp_m + m + \lambda
$. Since (\ref{B}) holds for all $y \in K$, we have $|Fp_m| \le \|p_m\| \le
3\varepsilon$.
Hence
$|k_m - m| \le |Fp_m| + \lambda \le 6\varepsilon,$
and thus
\begin{equation}\label{D}
k_m/m \to 1
\end{equation}
as $m \to\iy$.

Set $b_m = z + w_m/m$. Then $\| b_m \| \ge Fb_m = 1$. Since $\psi$ is an
$\varepsilon$-near\-isometry, we have
\begin{equation}\label{DD}
m - \varepsilon \le \| k_mz + w_m\| \le m + \varepsilon.
\end{equation}
Dividing by $m$ and letting $m \to \iy$ yields $\| b_m\| \to 1$ by 
(\ref{D}). Since $Fw_m =
0$, we have $\| z + w_m/2m\| \ge F(z + w_m/2m) = Fz = 1$. The uniform 
convexity of $X_2$ and
Lemma
\ref{RRR} imply that
\begin{equation}\label{E}
\| w_m\|/m = \| z-b_m\| \le \| b_m\| \gamma(1-1/\| b_m\|) \to 0,
\end{equation}
where $\gamma$ is as in (\ref{gamma}).

Write
$$
c = \psi (y) - P\psi (y),\quad v_m = \psi (y+mz),\quad P\psi (y) = \beta z 
+ u,
$$
where $\beta \in {\bf R} $ and $Fu=0$. Then $v_m = (\beta + k_m)z + w_m + u 
+ c$, and
hence $v_m/m \to z$ by (\ref{D}) and (\ref{E}). Since $X_2$ is uniformly 
smooth, this
implies that
$$
j(v_m) =j(v_m/m) \to j(z) = F
$$
by \cite[A.5, p. 411]{BeL}. Moreover, by (\ref{DD}) we have
$$
\| v_m\| =j(v_m)v_m \le j(v_m)(\beta z + u + c) + \| k_mz + w_m\| \le 
j(v_m)(\beta z + u
+ c) + m + \varepsilon.
$$
Since $Fc = FPc = 0$, this yields
\begin{equation}\label{F}
\limsup_{m \to\iy} (\| v_m\| - m) \le F(\beta z + u + c) + \varepsilon = 
\beta +
\varepsilon.
\end{equation}

Since $\psi$ is an $\varepsilon$-near\-isometry, we have
$\| y+mz\| - \| \psi (y+mz)\| \le \varepsilon,$
which can be rewritten as
$$
\| (\beta + \lambda + m)z + u\| - \| v_m\| \le \varepsilon.
$$
As $m \to\iy$, this implies by (\ref{A}) and (\ref{F}) that $\lambda - 
\varepsilon \le
\varepsilon$, and the theorem is proved.
\begin{rems} \rm
1. In particular, Theorem \ref{FigielLp} holds for
$\varepsilon$-near\-isometries $f: l^p
\to l^p$, $1 < p < \iy$.

2. The proof of Theorem \ref{FigielLp} is valid if
$X_1$ and
$X_2$ are Banach spaces and if (a) $X_2$ is uniformly convex and uniformly 
smooth, (b)
there is a linear projection of norm one of $X_2$ onto $\varphi (X_1)$, where 
$\varphi
(x) =
\lim_{s\to\iy} f(sx)/s$ exists by uniform convexity. The condition (b) is 
quite strong.
However, it is satisfied if $X_2$ is a Hilbert space, and we get the 
following result,
whose direct proof is given in \cite[5.3]{Vai}.
\end{rems}
\begin{izrek}
Let $H$ be a Hilbert space and $X$ any Banach space. Let $f: X \to H$ be an
$\varepsilon$-near\-isometry with $f(0) = 0$. Then there is a continuous 
linear
operator
$T:H \to X$ with $\| T \| = 1$ such that
$$
\| Tf(x)-x \| \le 2\varepsilon, \ \ \ x\in X.
$$
\end{izrek}

%\newpage

\section{Nearsurjective near\-isometries between general Banach spaces}
We start this section by recalling the following useful result of Dilworth 
\cite[Prop.
2]{Dil}.
\begin{lema}\label{di}
Suppose that $f: X \to Y$ is a near\-isometry between Banach spaces with 
$f(0)=0$ and that
$f$ maps
$X$
$\delta$-onto some closed subspace $Y_1$ of $Y$. Then there exist a number 
$L > 0$ and a
bijective linear iso\-metry $U:X \to Y_1$ such that $\| f(x) - Ux\| \le L$ 
for all $x \in
X$.
\end{lema}

Our first goal in this section is to improve the estimate (\ref{SD}) to the
sharp one.
\begin{izrek}\label{ISD}
Let $\varepsilon ,\delta \ge 0$ and let $f: X \to Y$ be an 
$\varepsilon$-near\-isometry
between Banach spaces satisfying $f(0) = 0$. Assume further that $f$ maps 
$X$ $\delta$-onto
some closed subspace $Y_1 \subset Y$. Then
there exists a bijective linear iso\-metry $U:X\to Y_1$ such that
\begin{equation}\label{2+4}
\| f(x) - Ux \| \le 2\varepsilon + 2\delta , \ \ \ x\in X.
\end{equation}

This estimate is sharp: For each pair $\varepsilon,\delta \ge 0$, there are Banach
spaces $X,Y$ and an $\varepsilon$-near\-isometry $f: X \to Y$ mapping $X$
$\delta$-onto a closed subspace $Y_1 \subset Y$ such that for each linear isometry
$U: X \to Y_1$, we have $\| f(x) - Ux\| \ge 2\varepsilon + 2\delta$ for some $x \in
X$.
\end{izrek}
{\bf Proof.}
As $f$ maps $X$
$\delta$-onto $Y_1$, we can find a map $g : X \to Y_1$
satisfying $g(0)=0$
such that
$\| f(x) - g(x) \| \le \delta$ for all $x\in X$. Then $g
:X\to Y_1$ is a $\delta_1$-near\-surjective
$\varepsilon_1$-near\-isometry with $\delta_1 = 2\delta$, 
$\varepsilon_1 = \varepsilon
+ 2\delta$.
By Lemma \ref{di} there exist a positive constant $L$
and a bijective linear iso\-metry $U:X\to Y_1$ such that
$
\| g(x) - Ux \| \le L
$ for all $ x\in X. $
Consequently,
$
\| f(x) - Ux \| \le L+\delta = K$ for all $x\in X$.

Replacing $f$ by $f U^{-1}: Y_1 \to Y$ we may assume that $X=Y_1 \subset Y$ 
and that $\|
f(x)-x\|
\le K$ for all $x \in X$. Let
$x
\in X$ and set $\lambda = \| g(x) - x\|$. It suffices to show that
\begin{equation}\label{a}
\lambda \le 2\varepsilon + \delta,
\end{equation}
because $\| f(x) - g(x)\| \le \delta$.

Choose a unit vector $u \in X$ with $x -g(x) = \lambda u$. For positive
integers $m$ we set $x_m = x + mu,\ z_m = f(x_m)-g(x),\ u_m = z_m/m$, and
\[
\alpha = \limsup_{m\to\iy}\, (\| f(x_m)\| - \| z_m\|).
\]
Observe that $|\alpha| \le \| g(x)\| < \iy$. Then $\| z_m - (m+\lambda)u\| 
= \| f(x_m)
- x_m\| \le K$ for all $m$. Dividing by $m$ we see that $u_m \to u$ as $m 
\to\iy$.

Let $0 < t < m$. Then $\| z_m\| > t\| u_m\|$, and thus $\| z_m\| = t\| 
u_m\| + \| f(x_m)
-(g(x)+tu_m)\|$. Consequently,
\[
\| f(x_m)\| - \| z_m\| \le \| g(x)+tu_m\| - t\| u_m\|.
\]
As $m\to\iy$, this yields $\alpha \le \| g(x)+tu\| - t$. Applying this with 
$t = \lambda +
n$, where $n$ is a positive integer, we get $\alpha \le \| x_n\| - \|
x_n-x\| -
\lambda$. Since $f$ is an $\varepsilon$-near\-isometry, this implies that
\[
\alpha \le \| f(x_n)\| - \| f(x_n) - f(x)\| + 2\varepsilon - \lambda \le \| 
f(x_n)\| - \|
z_n\| +
\delta + 2\varepsilon - \lambda.
\]
As $n \to\iy$, this gives \eqref{a}.

To prove the sharpness, let $X$ be the real axis, let $Y$ be the plane equipped
with the norm $\| (s,t)\| = |s| + |t|$, and define $f: X \to Y$ by setting $f(0) =
(0,0),\ f(\delta + \varepsilon) = (-\varepsilon,\delta)$, and
\begin{equation*}
f(t) = 
\begin{cases}
(t - \varepsilon,0)&  {\rm if}\ \ t <  0, \\
(-\varepsilon,t)& {\rm if}\ \ 0 < t\le \delta, \\
(t-\delta-\varepsilon, \delta)& {\rm if} \ \ t\ge \delta,\ t \ne
\delta+\varepsilon. 
\end{cases} 
\end{equation*}
Then  $f$ is an $\varepsilon$-near\-isometry with $f(0) = (0,0)$ mapping $X$
$\delta$-onto $Y_1 = \{  (s,t): t=0\}$.
There are only two linear isometries $U:X \to Y_1$, namely
$Ut = (t,0)$ and
$Ut = (-t,0)$. In the latter case $f-U$ is unbounded, and in the former case we
have
 $\| f(\delta+\varepsilon) - U(\delta+\varepsilon)	\| =
2\delta + 2
\varepsilon$.\\

In the case of Hilbert spaces we can improve the estimate of Theorem
\ref{ISD} as follows:

\begin{izrek}\label{Hilbert}
If the space $Y$ in Theorem \ref{ISD} is a Hilbert space, then the bound
$2\varepsilon + 2\delta$ in
\eqref{2+4}  can be
replaced by $2\varepsilon + \delta$ but not by $c\varepsilon + c'\delta$ for any
$c < 2$ or $c' < 1$.
\end{izrek}
{\bf Proof.}
We may again assume that $X = Y_1$ and that $\| f(x)-x\| \le K$ for all $x 
\in X$. Then
$f(sx)/s \to x$ as $s \to \iy$. Let $P: Y \to X$ be the orthogonal 
projection. Fix $x
\in X$, set $\alpha = \| Pf(x) - x\|$, and choose a unit vector $u \in X$ 
with $\alpha u
= Pf(x) - x$. By \cite[4.5]{Vai} we have
\[
|f(x) \cdot f(su) - x \cdot su| \le 2\varepsilon (\|x\| + \|su\| + 
\varepsilon),
\]
where $x \cdot y$ denotes the inner product. Dividing by $s$ and letting $s 
\to\iy$
yields $|f(x) \cdot u - x \cdot u| \le 2\varepsilon$. Since $f(x) \cdot u = 
Pf(x) \cdot
u$, this implies that $\alpha = \alpha u \cdot u \le 2\varepsilon$. Since 
$\| f(x) -
Pf(x)\| = {\rm dist}\, (f(x),X) \le \delta$, we obtain
\[
\| f(x)-x\| \le \alpha + \| f(x) - Pf(x)\| \le 2\varepsilon + \delta.
\]

The example  in \cite[p. 620]{OmS} or in \cite[p. 474]{Dil} shows that the
constant $2\varepsilon$ is optimal already in the case $\delta = 0$. To
prove that also the constant $\delta$ is optimal, let $X$ be the real line,
let $Y$ be the euclidean plane, let $\varepsilon,\delta > 0$, set $r =
\delta^2/(2\varepsilon)$, and define $f: X \to Y$ by 
\[ f(t) = 
\begin{cases} 0 &  {\rm if}\ \ t \le  0, \\  \delta t/r & {\rm
if}\ \ 0 \le t\le r, \\
\delta & {\rm if} \ \ t\ge r. 
\end{cases} 
\] Then $f$ is an $\varepsilon$-near\-isometry, and $f$ maps $X$
$\delta$-onto $Y_1 =
\{  (s,t): t = 0\}$. As in the proof of the previous theorem, it suffices to
consider the isometry $U: X \to Y_1$, $Ut = (t,0)$. For $t \ge r$ we have
$\|								 f(t) - Ut\| = \delta$. If the theorem holds with the bound
$c\varepsilon + c' \delta$,		we obtain $\delta \le c\varepsilon +
c'\delta$, which yields $c' \ge 1$ as $\varepsilon \to 0$.

\begin{rem} \rm
Although the bound $2\varepsilon + \delta$ in \ref{Hilbert} is the best
possible linear bound, the proof also gives the better nonlinear bound
$\sqrt{4\varepsilon^2 + \delta^2}$.
\end{rem}

We finally give the main result of the paper.
\begin{izrek}\label{IDT}
Suppose that $f: X \to Y$ is a
near\-surjective $\varepsilon$-near\-isometry
between Banach spaces satisfying $f(0) = 0$. Then
there exists a bijective linear iso\-metry $U:X\to Y$ such that
\begin{equation}\label{twotwo}
\| f(x) - Ux \| \le 2\varepsilon,\ x \in X.
\end{equation}
Hence $f$ is $2\varepsilon$-near\-surjective. The constant $2$ is the best 
possible.
\end{izrek}
{\bf Proof.} By Lemma \ref{di} there exist a positive constant $L$
and a bijective linear iso\-metry $U:X\to Y$ such that
$
\| f(x) - Ux \| \le L
$ for all $ x\in X. $
Replacing $f$ by $f U^{-1}: Y \to Y$ we may assume that $X= Y$ and that $\|
f(x)-x\|
\le L$ for all $x \in X$.

We can now proceed as in the proof of 
Theorem \ref{ISD}, but the situation is 
simpler, because
$g(x)$ is replaced by
$f(x)$. For $\lambda = \| f(x)-x\|$ we obtain the inequality 
$\alpha \le \|
f(x_n)\| -
\| z_n\| + 2\varepsilon - \lambda$, which gives $\lambda \le 2\varepsilon$ 
as $n \to\iy$. 

The sharpness of the bound is well known already in the case $\delta = 0$;
see \cite[p. 620]{OmS} or \cite[p. 474]{Dil}.
\begin{rem}\rm
Theorem \ref{IDT} was proved by Dilworth \cite[Th. 2]{Dil} under the 
additional
condition that the norm of
$X$ is Fr\' echet differentiable at each point of  a dense set.
\end{rem}
\baselineskip16pt
{\small
% -----------------------------------------------------------------------

}

\noindent Peter \v Semrl \hfill Jussi V\" ais\" al\" a\\
Department of Mathematics  \hfill Matematiikan laitos\\
University of Ljubljana \hfill Helsingin yliopisto\\
Jadranska 19 \hfill PL 4, Yliopistonkatu 5\\
SI-1000 Ljubljana, Slovenia  \hfill  FIN-00014 Helsinki, Finland\\
\texttt{peter.semrl@fmf.uni-lj.si  \hfill  jvaisala@cc.helsinki.fi}\\

\end{document}